\documentclass[a4 paper, 12pt, reqno]{amsart}
\usepackage{latexsym}
\usepackage[T1]{fontenc}
\usepackage[english]{babel}
\usepackage{amssymb,amsmath,amsthm,amsfonts}
\usepackage[latin1]{inputenc}

\usepackage{url}
\usepackage{a4wide}
\usepackage{enumerate}
\usepackage{changebar}
\usepackage{booktabs}
\usepackage{paralist}
\usepackage{comment}
\usepackage{mathrsfs}

\usepackage{graphicx}
\usepackage{color}
\usepackage{cite}
\theoremstyle{plain}
\newtheorem{theorem}{Theorem}[section]
\newtheorem{prop}[theorem]{Proposition}
\newtheorem{lemma}[theorem]{Lemma}
\newtheorem{corollary}[theorem]{Corollary}
\theoremstyle{definition}
\newtheorem{definition}[theorem]{Definition}
\newtheorem{remark}[theorem]{Remark}

\def\R{{\mathbb R}}

\def\N{{\mathbb N}}

\def\D{{\mathcal{D}}}

\def\OU{\mathcal{L}}
\def\e{\mathrm{e}}
\def\L{\mathscr{L}}
\def\B{\mathrm{B}}
\def\d{\mathrm{d}}

\numberwithin{equation}{section}

\title[Non-autonomous Ornstein-Uhlenbeck equations ]{\textbf{Non-autonomous Ornstein-Uhlenbeck equations in exterior domains}}
\author{Tobias Hansel}
\address{ Technische Universität Darmstadt\\Department of Mathematics\\ 64289 Darmstadt,
Germany} \email{hansel@mathematik.tu-darmstadt.de}
\author{Abdelaziz Rhandi}
\address{Dipartimento di Ingegneria dell'Informazione e Matematica Applicata \\ Universit\`a degli Studi di Salerno\\ Via Ponte Don Melillo, \\ 84084 Fisciano
(Sa)\\Italy}
\email{rhandi@diima.unisa.it}
\begin{document}

%%%%%%%%%%%%%%%%%%%%%%%%%%%%%%%%%%%%%%%%%%%%%%%%%%%%%%%%%%%%%%%%%%%%

\maketitle %\thispagestyle{fancy}
%\tableofcontents
%\today

%%%%%%%%%%%
%
\begin{abstract}
In this paper, we consider non-autonomous Ornstein-Uhlenbeck  operators
in smooth exterior domains $\Omega\subset \R^d$ subject to  Dirichlet
boundary conditions. Under suitable assumptions on the coefficients,
the solution of the corresponding non-autonomous parabolic Cauchy problem is governed by an evolution
system $\{P_\Omega(t,s)\}_{0\le s\le t}$ on
$L^p(\Omega)$ for $1< p < \infty$.
Furthermore, $L^p$-estimates for spatial
derivatives and $L^p$-$L^q$ smoothing properties of $P_\Omega(t,s),\,0\le s\le t,$ are obtained.
\end{abstract}

\section{Introduction}
In recent years, parabolic equations with unbounded and time-independent coefficients were investigated intensively in various function 
spaces over the whole space $\R^d$ or exterior domains; we refer e.g. to \cite{DaPrato/Lunardi:1995,Geissert/Heck/Hieber:
2006,Geissert/etal:2005, Hieber/etal:2007,Hieber-Sawada} and the monograph \cite{Bertoldi/Lorenzi}. However, it is also interesting to 
consider parabolic 
equations with unbounded coefficients in the non-autonomous case. In particular, analytically  there is a great interest in the prototype 
situation of time-dependent Ornstein-Uhlenbeck operators in exterior domains, as operators of this type arise e.g. in the study of the 
Navier-Stokes flow in the exterior of a rotating obstacle; see e.g. \cite{Hansel:2009,Hishida:2001}.

Therefore, in this paper we consider non-autonomous Cauchy problems with Dirichlet boundary condition of the type
\begin{equation}\label{eq:FCP}
\left\{
\begin{array}{rclll}
u_t(t,x)- \OU_\Omega(t)u(t,x)&=&0, &t\in (s,\infty),\; x \in \Omega, \\[0.2cm]
u(t,x)&=& 0, & t\in (s,\infty),\; x \in \partial\Omega, \\[0.2cm]
u(s,x) & = & f(x), &  x\in\Omega,
\end{array}\right.
\end{equation}
where $s\geq 0$ is fixed, $\Omega\subset \R^d$ is a domain and $\{\OU_\Omega(t)\}_{t\ge 0}$ is a family of time-dependent Ornstein-Uhlenbeck operators formally defined by
\begin{equation}
\OU_\Omega(t) \varphi(x)= \frac 1 2 \mathrm{Tr}\left(Q(t)Q^*(t)\mathrm D_x^2 \varphi(x)\right) +  \langle M(t)x + c(t), \mathrm D_x\varphi(x) \rangle, \quad x \in \Omega, \quad t\ge 0.
\end{equation}
Throughout the paper we assume that $Q,\,M\in C^\alpha_{loc}(\R_+, \R^{d\times d}),\,c\in C^\alpha_{loc}(\R_+, \R^{d})$ for some $\alpha \in (0,1)$ and there is $\mu >0$ such that
$$|Q(t)x|\ge \mu |x|,\quad t\ge 0,\,x\in \R^d.$$ The above assumption guaranties that the operators $\OU_\Omega(t)$ are uniformly elliptic.

The main purpose of this paper is to consider problem \eqref{eq:FCP} in the $L^p$-setting for the case of smooth exterior domains
$\Omega$. However, in the course of this paper we also consider the situation where $\Omega$ is $\R^d$ and a smooth bounded domain.

In the following the $L^p$-realization of $\OU_\Omega(t)$ will be denoted by $L_\Omega(t)$ with an appropriate domain $\D(L_
\Omega(t))\subset L^p(\Omega)$, specified later. Then we can rewrite equation (\ref{eq:FCP}) as an abstract non-autonomous Cauchy
problem
\begin{equation}\label{eq:nACP}
\mathrm{(nACP)}\qquad\left\{
\begin{array}{rclll}
u'(t)& = & L_{\Omega}(t)u(t), & 0\leq s< t , \\[0.2cm]
u(s) & = & f,&
\end{array}\right.
\end{equation}
where $f  \in L^p(\Omega)$.
\begin{definition}
A continuous function $u:[s,\infty) \to L^p(\Omega)$ is called a \emph{(classical) solution} of (nACP) if $u\in C^1((s,\infty), L^p(\Omega))$, $u(s)=f$, and $u'(t)=L_\Omega(t)u(t)$ for $0\leq s< t$.
\end{definition}
\begin{definition}[Well-posedness]
We say that
the Cauchy problem (nACP) is \emph{well-posed} (on \emph{regularity spaces} $\{Y_s\}_{s\geq 0}$) if the following statements are true.
\begin{itemize}
\item[(i)] {\bf (Existence and uniqueness)} There are dense subspaces $Y_s \subset \D(L_\Omega(s))$ of $L^p(\Omega)$
such that for $f\in Y_s$ there is a unique solution $t
\mapsto u(t;s,f) \in Y_t$ of (nACP).
\item[(ii)] {\bf (Continuous dependence)} The solution depends continuously on the data; i.e., for $s_n
\rightarrow s$ and $Y_{s_n} \ni f_{n} \rightarrow f \in Y_s$, we have
$\tilde u(t;s_n,f_{n})\rightarrow \tilde u(t;s,f)$ uniformly for $t$
in compact subsets of $[0,\infty)$, where we set $\tilde
u(t;s,f):=u(t;s,f)$ for $t\geq s$ and $\tilde u(t;s,f):=f$ for
$t<s$.
\end{itemize}
\end{definition}
In order to discuss well-posedness of (nACP) we introduce the concept of strongly continuous evolution systems.
\begin{definition}[Evolution system]\label{def:evolution_system}
A two parameter family of linear, bounded operators $\{P_{\Omega}(t,s)\}_{0\leq s \leq t }$ on $L^p(\Omega)$ is called a \emph{(strongly continuous) evolution system}
if
\begin{itemize}
\item[(i)] $P_{\Omega}(s,s)=\mathrm{Id}$ \hspace{0.4cm} and \hspace{0.4cm} $P_{\Omega}(t,s) = P_{\Omega}(t,r)P_{\Omega}(r,s)$ \hspace{0.6cm} for $0 \leq s \leq r \leq t$,
\item[(ii)] for each $f\in L^p(\Omega)$, \hspace{0.4cm} $(t,s) \mapsto
P_{\Omega}(t,s)f$ \hspace{0.6cm} is continuous on $0 \leq s \leq t$.
\end{itemize}
We say $\{P_{\Omega}(t,s)\}_{0\leq s \leq t }$ solves the Cauchy problem (nACP) (on
spaces $\{Y_s\}_{s\geq 0}$) if there are dense subspaces $Y_s$ of $L^p(\Omega)$
such that $P_{\Omega}(t,s)Y_s \subset Y_t \subset \D(L_\Omega(t))$ for $0\leq s\leq t$  and
the function $u(t):= P_{\Omega}(t,s)f$ is a solution of (nACP) for $f\in Y_s$.
\end{definition}
It is well-known that the Cauchy problem (nACP) is well-posed on $\{Y_s\}_{s\geq 0}$ if and only if
there is an evolution system solving (nACP) on $\{Y_s\}_{s\geq 0}$ (see e.g. \cite[Sect. 3.2]{Nickel:1996}).

The main result of this paper (see Theorem \ref{prop:evolution_system_exterior}) is to show that for smooth exterior domains $\Omega
\subset\R^d$ problem (nACP) is solved by a strongly continuous evolution system $\{P_{\Omega}(t,s)\}_{0\leq s \leq t }$ on 
$L^p(\Omega)$ and thus, is well-posed. Since in unbounded domains the operators $\OU_\Omega(t)$ have unbounded drift coefficients, 
the present situation does not fit into the well-studied framework of evolution systems of parabolic type (see e.g. the monograph by 
Lunardi \cite[Chapter 6]{Lunardi:1995} or the fundamental papers by Tanabe \cite{Tanabe:1959, Tanabe:1960a, Tanabe:1960b} and 
Acquistapace, Terreni \cite{Acquistapace:1984,Acquistapace/Terreni:1986,Acquistapace/Terreni:1987}). Therefore the well-posedness 
of (nACP) and regularity properties of the solution do not follow from abstract arguments. Here lies the major difficulty. In order to prove 
our result we proceed as follows: In Section 2 we consider (nACP) in the case that $\Omega$ is the whole space $\R^d$ or a smooth 
bounded domain. For the whole space case we use a representation formula for the evolution system as done in \cite{DaPrato/Lunardi:
2007, Geissert/Lunardi:2008}. In the case of bounded domains we can apply the standard results for non-autonomous Cauchy problems 
of parabolic type. These auxiliary results are then applied in Section 3 to construct an evolution system $\{P_{\Omega}(t,s)\}_{0\leq s \leq 
t } $ on $L^p(\Omega)$ for smooth exterior domains $\Omega\subset \R^d$, by some cut-off techniques. Moreover, our method allows us 
to prove $L^p$-$L^q$ estimates and estimates for spatial derivatives of $\{P_{\Omega}(t,s)\}_{0\leq s \leq t }$.

\subsection*{Notations}
The euclidian norm of $x\in \R^d$ will be denoted by $|x|$. By $B(R)$ we denote the open ball in $\R^d$ with centre
at the origin and radius $R$. For $T>0$ we use the notations:
\begin{eqnarray*}
\Lambda_T &:=& \{(t,s): 0\le s\le t\le T\}\\
\widetilde{\Lambda}_T &:=&  \{(t,s): 0\le s< t\le T\}\\
\Lambda &:=& \{(t,s): 0\le s\le t \}\\
\widetilde{\Lambda} &:=&  \{(t,s): 0\le s< t \}.
\end{eqnarray*}
If $u: \Omega \to \R$, where $\Omega \subseteq \R^d$ is a domain,
we use the following notation:
\begin{eqnarray*}
\mathrm D_i u &=& \frac{\partial u}{\partial x_i}, \
\mathrm D_{ij} u=\mathrm D_i \mathrm D_j u,
\\
\mathrm D_x u&=&(\mathrm D_1 u, \dots, \mathrm D_d u), \ \mathrm D^2_x u=(\mathrm D_{ij} u).
\end{eqnarray*}
Let us come to notation for function spaces. For $1\le p< \infty,\,j\in \N,\,W^{j,p}(\Omega)$ denotes the
classical Sobolev space of all $L^p(\Omega)$--functions having weak derivatives in
$L^p(\Omega)$ up to the order $j$.  Its usual norm is denoted by $\|\cdot
\|_{j,p}$ and by $\|\cdot \|_p$ when $j=0$. By $W_0^{1,p}(\R^d)$ we denote the closure of the space of test functions $C_c^\infty(\R^d)$ with respect to the norm of $W^{1,p}(\R^d)$. For $0<\alpha < 1$ we denote by $C^\alpha_{loc}(\R_+, \R^{d\times d})$ the space of all $\alpha$-H\"older continuous functions in $[0,T]$ for all $T>0$. The space of all bounded continuous functions $u:\Omega\to \R$ is denoted by $C_b(\Omega)$. For $k\in \N$, $C_b^k(\Omega)$ is the subspace of $C_b(\Omega)$ consisting of all functions which are differentiable up to the order $k$ in $\Omega$ such that the derivatives are bounded. Finally, we denote by $C^{1,2}(I\times \Omega)$ the space of all functions $u: I\times \Omega \to \R$ which are continuously differentiable with respect to $t\in I$ and $C^2$ with respect to the space variable $x\in \Omega$, where $I\subseteq [0,\infty)$ is an interval.
\section{Auxiliary results: whole space and bounded domains}
In this section we prove some auxiliary results concerning the evolution systems in the case of the whole space $\R^d$ and smooth bounded domains. These results are needed in Section 3 for the construction of the evolution system in the case of exterior domains.
\subsection{The evolution system in the whole space}
The realizations of $\{\OU_{\R^d}(t)\}_{t\ge 0}$ are defined by
\begin{equation}
\begin{array}{lcl}
\D(L_{\R^d}(t))&:=& \{u \in W^{2,p}(\R^d): \langle M(t)x, \mathrm D_xu(x)
\rangle \in L^p(\R^d) \},\\[0.2cm]
L_{\R^d}(t)u &:=& \OU_{\R^d}(t) u.
\end{array}
\end{equation}
Here the domain of $L_{\Omega}(t)$ depends on the time parameter $t$. However, note that the subspace
\begin{equation*}
Y_{\R^d}:=\{u \in W^{2,p}(\R^d) : |x|\cdot \mathrm D_j u(x) \in L^p(\R^d) \; \hbox{\ for all} \; j=1,\ldots, d\}
\end{equation*}
is contained in $\D(L_\Omega(t))$ for all $t\ge 0$ and is dense in $L^p(\R^d)$. The space $Y_{\R^d}$ will serve as a regularity space in order to discuss well-posedness of (nACP).

It follows directly from \cite{Metafune/etal:2002} (see also \cite{Metafune:2001}) that in the autonomous case (i.e. for fixed $s\ge 0$) the operator $(L_{\R^d}(s),\D(L_{\R^d}(s))$ generates a strongly continuous semigroup, which is however not analytic. Second order elliptic operators in $\R^d$ with more general unbounded and time-independent coefficients were considered e.g. in  \cite{Pruess/etal:2006},   \cite{Hieber/etal:2009}.

In the following we denote by $\{U(t,s)\}_{t,s\ge 0}$ the evolution system in $\R^d$ that satisfies
\begin{equation*}
\left\{\begin{array}{lcl}
\frac{\partial}{\partial t} U(t,s) & = & -M(t)U(t,s),\\[0.2cm]
U(s,s) & = & \mathrm{Id}.
\end{array}\right.
\end{equation*}
The existence of $\{U(t,s)\}_{t,s\ge 0}$ follows directly from the Picard-Lindel\"of theorem. Now for $f\in L^p(\R^d)$ and $s\geq 0$ we set $P_{\R^d}(s,s)=\mathrm{Id}$ and for $(t,s)\in \widetilde{\Lambda}$ we define
\begin{equation}\label{eq:evol_sys_whole}
P_{\R^d}(t,s)f(x) = (k(t,s,\cdot)\ast f)(U(s,t)x+g(t,s)), \qquad x\in \R^d,
\end{equation}
where
\begin{equation}
k(t,s,x):= \frac{1}{(2\pi)^{\frac d 2} (\det Q_{t,s})^{\frac 1 2}} \ \mathrm{e}^{-\frac 1 2 \langle Q_{t,s}^{-1}x,x\rangle}, \qquad x\in \R^d,
\end{equation}
\begin{equation}
g(t,s) = \int_s^t U(s,r)c(r) \d r \qquad \mbox{and} \qquad Q_{t,s}=\int_s^t U(s,r)Q(r)Q^*(r)U^*(s,r)\d r.
\end{equation}
As in \cite[Proposition 2.1]{DaPrato/Lunardi:2007} (see also \cite[Proposition 2.1]{Hansel:2009}) it can be shown that for initial value $f\in C^2_b(\R^d)$, the function $u(t,x):= P_{\R^d}(t,s)f(x)$ is a classical solution to
\begin{equation}\label{eq2.5}
\left\{
\begin{array}{rclll}
u_t(t,x)- \OU_{\R^d}(t)u(t,x)&=&0, &(t,s)\in \widetilde{\Lambda},\; x \in \R^d, \\[0.2cm]
u(s,x) & = & f(x), &  x\in\R^d ,
\end{array}\right.
\end{equation}
i.e. $u\in C^{1,2}((s,\infty)\times \Omega)$ and $u$ solves (\ref{eq2.5}).
Further, the two parameter family of operators $\{P_{\R^d}(t,s)\}_{(t,s)\in \Lambda}$ is a strongly continuous evolution system on $L^p(\R^d)$.
\begin{prop}\label{prop:evolution_whole}
Let $1<p<\infty$. Then the family of operators $\{P_{\R^d}(t,s)\}_{(t,s)\in \Lambda}$ defined in \eqref{eq:evol_sys_whole} is a strongly continuous evolution system on $L^p(\R^d)$ with the following properties.
\begin{itemize}
\item[(a)] For $(t,s)\in \Lambda$, the operator $P_{\R^d}(t,s)$ maps $Y_{\R^d}$ into $Y_{\R^d}$.
\item[(b)] For every $f\in Y_{\R^d}$ and every $s\in[0,\infty)$, the map $t\mapsto P_{\R^d}(t,s)f$ is differentiable in $(s,\infty)$ and
\begin{equation}\label{eq:derivative_t}
\frac{\partial}{\partial t} P_{\R^d}(t,s)f = L_{\R^d}(t)P_{\R^d}(t,s)f .
\end{equation}
\item[(c)] For every $f\in Y_{\R^d}$ and $t\in (0,\infty)$, the map $s\mapsto P_{\R^d}(t,s)f$ is differentiable in $[0,t)$
and
\begin{equation}\label{eq:derivative_s}
\frac{\partial}{\partial s} P_{\R^d}(t,s)f = -P_{\R^d}(t,s)L_{\R^d}(s)f .
\end{equation}
\end{itemize}
\end{prop}
\begin{proof}
In \cite[Proposition 2.4]{Geissert/Lunardi:2008} it was shown
that the law of evolution  (property (i) of Definition \ref{def:evolution_system})
holds for every $f \in C_c^{\infty}(\R^d)$. Since $C_c^{\infty}(\R^d)$ is dense in $L^p(\R^d)$ the law of evolution holds even for all $f\in L^p(\R^d)$. The strong continuity of the map $\Lambda \ni (t,s)\mapsto P_{\R^d}(t,s)$ can be shown as in \cite[Proposition 2.3]{Hansel:2009}. Equalities \eqref{eq:derivative_t} and \eqref{eq:derivative_s} follow by differentiating the kernel $k(t,s,x)$ with respect to $t$ and $s$, respectively.

Let us now show that the evolution system $\{P_{\R^d}(t,s)\}_{(t,s)\in \Lambda }$ leaves the regularity space $Y_{\R^d}$ invariant. Since $k(t,s,\cdot) \in C^\infty(\R^d)$ it follows that $P_{\R^d}(t,s)f \in C^\infty(\R^d)$ for all $f\in L^p(\R^d)$ and $(t,s)\in \widetilde{\Lambda}$. Moreover, we note that
\begin{equation*}
\mathrm D_x P_{\R^d}(t,s)f = U^*(s,t) \left( k(t,s,\cdot)\ast \mathrm D_x f\right)(U(s,t)x+g(t,s))
\end{equation*}
holds for all $f\in W^{1,p}(\R^d)$. Thus, it suffices to show that for all $j=1,\ldots,d$ we have $|x| \cdot \left( k(t,s,\cdot)\ast \mathrm D_j f\right)(x)\in L^p(\R^d)$. So let $h\in L^q(\R^d)$ with $\frac 1 p + \frac 1 q = 1$. Then we obtain
\begin{align*}
& \int_{\R^d}\big| \left(|x| \cdot \left( k(t,s,\cdot)\ast \mathrm D_j f\right)(x)\right) h(x) \big|\d x\\
& \qquad \qquad\leq C \int_{\R^d} |x| |h(x)| \int_{\R^d} |\mathrm D_j f(x-y) \e ^{-\frac 1 2 \langle Q_{t,s}^{-1}y,y \rangle}|\d y \;\d x\\
& \qquad \qquad \leq C \left[\int_{\R^d} \e ^{-\frac 1 2 \langle Q_{t,s}^{-1}y,y \rangle} \int_{\R^d} \big|\left(|x-y|\cdot \mathrm D_j f(x-y)\right) h(x)\big|\d x \;\d y +\right. \\
& \qquad \qquad \qquad \qquad \left.\int_{\R^d}|y| \e ^{-\frac 1 2 \langle Q_{t,s}^{-1}y,y \rangle} \int_{\R^d}|\mathrm D_j f(x-y)||h(x)|\d x \;
\d y \right] \\
& \qquad \qquad \le C\left[\,\||x|\mathrm D_j f\|_p\|h\|_q+\|\mathrm D_j f\|_p\|h\|_q\,\right].
\end{align*}
Here the constant $C$ may change from line to line.
 Thus
\begin{equation*}
 \int_{\R^d}\big| \left(|x| \cdot \left( k(t,s,\cdot)\ast \mathrm D_j f\right)(x)\right) h(x) \big|\d x < \infty
\end{equation*}
holds for all $h \in L^q(\R^d)$ and this proves the assertion.
\end{proof}
As a consequence of Proposition \ref{prop:evolution_whole}, Cauchy problem (nACP) is well-posed in the case of $\R^d$ with regularity space $Y_{\R^d}$. Now we prove $L^p$-$L^q$ estimates and estimates for higher order spatial derivatives of $\{P_{\R^d}(t,s)\}_{(t,s)\in \Lambda}$. For this purpose we need the following estimates for the matrices $Q_{t,s}$. For a proof we refer to \cite[Lemma
3.2]{Geissert/Lunardi:2008} and \cite[Lemma 2.4]{Hansel:2009}.
\begin{lemma}\label{lemma:Qts_Estimates}
Let $T>0$. Then there exists a constant $C:=C(T)>0$ such that
\begin{align}
	\begin{split}
\|Q_{t,s}^{-\frac 1 2}\| \leq C (t-s)^{-\frac 1 2}, \quad (t,s)\in \widetilde{\Lambda}_T,\\
(\det Q_{t,s})^{\frac 1 2} \geq
C(t-s)^{\frac d 2},  \quad (t,s)\in \Lambda_T.
\end{split}
\label{eq:estimateQ}
\end{align}
\end{lemma}
\begin{prop}\label{prop:Lp_Lq_estimates_whole}
Let $T>0$, $1<p\leq q <\infty$ and $\beta \in \N_0^d$ be a multi-index . Then there exists a constant $C:=C(T)>0$ such that for every $f\in L^p(\R^d)$
\begin{enumerate}
\item[(a)] $\|P_{\R^d}(t,s)f\|_{q} \leq C(t-s)^{-\frac{d}{2}\left(\frac 1 p - \frac 1 q\right)}\|f\|_{p}, \qquad (t,s)\in \widetilde{\Lambda}_T$,\vspace{0.2cm}
\item[(b)]  $\|\mathrm D_x^\beta P_{\R^d}(t,s)f\|_{p} \leq C(t-s)^{-\frac{|\beta|}{2}}\|f\|_{p}, \qquad (t,s)\in \widetilde{\Lambda}_T$.
\end{enumerate}
Moreover, $$\|P_{\R^d}(t,s)f\|_{k,p}\le C\|f\|_{k,p},\quad (t,s)\in \Lambda_T,$$
for all $f\in W^{k,p}(\R^d),\,k=1,\,2,$ and
$$\|P_{\R^d}(t,s)f\|_{2,p}\le C(t-s)^{-\frac{1}{2}}\|f\|_{1,p},\quad (t,s)\in \widetilde{\Lambda}_T,$$
for all $f\in W^{1,p}(\R^d)$.
\end{prop}
\begin{proof}
Let $T>0$. By the change of variables $\xi=U(s,t)x$ and by Young's inequality we obtain
\begin{equation*}
\|P_{\R^d}(t,s)f\|_{q}\leq |\det U(s,t)|^{\frac 1 q} \|k(t,s,\cdot)\|_{r} \|f\|_{p},\vspace{0.05cm}
\end{equation*}
where $1<r<\infty$ with $\frac 1 p + \frac 1 r = 1 + \frac 1 q$. Moreover, by the change of variables $y=Q_{t,s}^{1/2}z$ we obtain
\begin{equation*}
\|k(t,s,\cdot)\|^r_{r} = \frac{(\det Q_{t,s})^{\frac 1 2 (1-r)}}{(2\pi)^{\frac d 2 \cdot r}} \int_{\R^d}\e^{-\frac{r |z|^2}{2}}\d z \leq C (\det Q_{t,s})^{\frac 1 2 (1-r)}.\vspace{0.05cm}
\end{equation*}
Now Lemma \ref{lemma:Qts_Estimates} yields (a).

To prove (b) we first note that
\begin{equation*}
\left|\mathrm D_x^{\beta}P_{\R^d}(t,s)f(x)\right| \le \left|U^*(s,t)\right|^{|\beta|} \left|\left(\mathrm D_x^\beta k(t,s,\cdot)\ast f \right)(U(s,t)x+g(t,s))\right|
\end{equation*}
holds. Thus, we have to estimate the norm of $\mathrm D_x^\beta k(t,s,\cdot)$. Since
\begin{equation*}
\mathrm D_x k(t,s,x)= - k(t,s,x)\left(Q_{t,s}^{-1}x \right)^*
\end{equation*}
holds, we obtain by differentiating further
\begin{equation*}
|\mathrm D_x^{\beta} k(t,s,x)|\leq C k(t,s,x) |Q_{t,s}^{-1}x|^{|\beta|}
\end{equation*}
for some constant $C>0$. As above, by the change of variables $y=Q_{t,s}^{1/2}z$, we obtain
\begin{equation*}
\|\mathrm D_x^{\beta} k(t,s,\cdot) \|_{1} \leq \frac{\|Q_{t,s}^{-\frac{1}{2}}\|^{|\beta|}}{(2\pi)^{\frac d 2}} \int_{\R^d}|z|^{|\beta|}\e^{-\frac{|z|^2}{2}}\d z \leq C\|Q_{t,s}^{-\frac{1}{2}}\|^{|\beta|}.
\end{equation*}
Now Lemma \ref{lemma:Qts_Estimates} yields assertion (b). 
The last assertions follow by a direct computation.
\end{proof}
\begin{remark}
If $\{U(t,s)\}_{t,s\ge 0}$ is uniformly bounded, i.e. $\|U(t,s)\|\leq M$ for some constant $M>0$ and all $t,s\ge 0$, then the estimates in Lemma \ref{lemma:Qts_Estimates} and Proposition \ref{prop:Lp_Lq_estimates_whole} hold in $\Lambda$ and $\widetilde{\Lambda}$ respectively. In particular, in this case the evolution system $\{P(t,s)\}_{(t,s)\in\Lambda}$ is uniformly bounded.
\end{remark}
\subsection{The evolution system in bounded domains}
In this subsection we assume that $D\subset \R^d$ is a bounded domain with   $C^{1,1}$-boundary. For $t\ge 0$ we set
\begin{equation}
\begin{array}{lcl}
\D(L_{D}(t))&=:& \D(L_D):= W^{2,p}(D)\cap W^{1,p}_0(D),\\[0.2cm]
L_{D}(t)u &:=& \OU_D(t) u.
\end{array}
\end{equation}
Note that in this situation the domain is independent of the time parameter $t$, i.e. all the operators $L_{D}(t)$ are defined on the same domain $\D(L_D)$.
\begin{lemma}
Let $D \subset \R^d$ be a bounded domain with $C^{1,1}$-boundary and $1<p<\infty$.
\begin{enumerate}
\item[(a)] For fixed $s\in[0,\infty)$, the operator $(L_{D}(s),\D(L_D))$ generates an analytic semigroup on $L^p(D)$. 
\item[(b)] The map $t\mapsto L_{D}(t)$ belongs to $C^\alpha_{loc}(\R_+, \L(\D(L_D),L^p(D)))$.
\end{enumerate}
\end{lemma}
\begin{proof}
Assertion (a) follows from the classical theory of elliptic second order operators  in bounded domains (see also \cite[Lemma 2.4]{Geissert/etal:2005}). Assertion (b) follows from the assumptions on the coefficients of $L_{D}(\cdot)$.
\end{proof}
The following proposition now follows directly from the theory of evolution systems of parabolic type; see \cite[Chapter 6]{Lunardi:1995} and \cite[Sect. 2.3]{Grisvard}. See also \cite[Sect. 7]{Amann} for bounded domains of class $C^2$.
\begin{prop}\label{prop:evolution_bounded}
Let $D\subset \R^d$ be a bounded domain with $C^{1,1}$-boundary and $1<p<\infty$. Then there is a unique evolution system $\{P_{D}(t,s)\}_{(t,s)\in \Lambda}$ on $L^p(D)$ with the following properties.
\begin{enumerate}
\item[(a)] For $(t,s)\in \widetilde{\Lambda}$, the operator $P_{D}(t,s)$ maps $L^p(D)$ into $\D(L_D)$.
\item[(b)] The map $t \mapsto P_{D}(t,s)$ is differentiable in $(s,\infty)$ with values in $\L(L^p(D))$ and
\begin{equation}
\frac{\partial}{\partial t} P_{D}(t,s) = L_{D}(t)P_{D}(t,s).
\end{equation}
\item[(c)] For every $f \in \D(L_D)$ and $t\in (0,\infty)$, the map $s\mapsto P_{D}(t,s)f$ is differentiable in $[0,t)$ and
\begin{equation}
\frac{\partial}{\partial s} P_{D}(t,s)f = - P_{D}(t,s)L_{D}(s)f.
\end{equation}
\item[(d)] Let $T>0$. Then there exists a constant $C:=C(T)>0$ such that
\begin{equation}
\|P_{D}(t,s)f\|_{p}\leq C\|f\|_{p},
\end{equation}
and
\begin{equation}
\|P_{D}(t,s)f\|_{2,p}\leq C(t-s)^{-1}\|f\|_{p}.
\end{equation}
for all $f\in L^p(D)$ and all $(t,s)\in \widetilde{\Lambda}_T$.
\end{enumerate}
\end{prop}

The following estimates follow directly from the proposition above and simple interpolation.
\begin{corollary}\label{prop:Lp_Lq_estimates_bounded}
Let $T>0$, $1<p<\infty$ and $p\leq q < \infty$. Then there exists a constant $C:=C(T)>0$ such that for every $f\in L^p(D)$
\begin{enumerate}
\item[(a)] $\|P_{D}(t,s)f\|_{q} \leq C(t-s)^{-\frac{d}{2}\left(\frac 1 p - \frac 1 q\right)}\|f\|_{p}, \qquad (t,s)\in \widetilde{\Lambda}_T$,\vspace{0.2cm}
\item[(b)]  $\|\mathrm D_x P_{D}(t,s)f\|_{p} \leq C(t-s)^{-\frac 1 2}\|f\|_{p}, \qquad (t,s)\in \widetilde{\Lambda}_T$.
\end{enumerate}
Moreover, $$\|P_{D}(t,s)f\|_{k,p}\le C\|f\|_{k,p},\quad (t,s)\in \Lambda_T,$$
for all $f\in W^{k,p}(D),\,k=1,\,2,$ and
$$\|P_{D}(t,s)f\|_{2,p}\le C(t-s)^{-\frac{1}{2}}\|f\|_{1,p},\quad (t,s)\in \widetilde{\Lambda}_T,$$
for all $f\in W^{1,p}(D)$.
\end{corollary}
\begin{proof}
Let us start with the case $q\geq p \geq d/2$. Then, by the Gagliardo-Nierenberg inequality (cf. \cite[Theorem 3.3]{Tanabe:1997}) and Proposition \ref{prop:evolution_bounded} (d), we immediately obtain
\begin{align*}
\|P_{D}(t,s)f\|_{q} &\leq C \|\mathrm D_x^2 P_{D}(t,s)f\|_{p}^a \|P_{D}(t,s)f\|_{p}^{1-a} \leq C(t-s)^{-a}\|f\|_{p},\,(t,s)\in \widetilde{\Lambda}_T,
\end{align*}
where $a= \frac d 2 \left(\frac 1 p - \frac 1 q \right)$. The case $1<p<\frac d 2$ follows by iteration. Assertion (b) is also proved by the Gagliardo-Nierenberg inequality. By setting $a=\frac 1 2$ and $p=q$ we obtain
\begin{align*}
\|\mathrm D_x P_{D}(t,s)f\|_{p} &\leq C \|\mathrm D_x^2 P_{D}(t,s)f\|_{p}^{\frac 1 2 } \|P_{D}(t,s)f\|_{p}^{\frac 1 2 } \leq C(t-s)^{-\frac 1 2}\|f\|_{p},\,(t,s)\in \widetilde{\Lambda}_T.
\end{align*}
For the last assertions we refer, for example, to \cite[Corollary 6.1.8]{Lunardi:1995}.
\end{proof}
\section{The evolution system in exterior domains}
In this section we come to the main part of this paper. In the sequel we always assume that $\Omega\subset \R^d$ is an exterior domain with $C^{1,1}$-boundary, i.e., $\Omega = \R^d \setminus K$, where $K\subset \R^d$ is a compact set with $C^{1,1}$-boundary. For $t\ge 0$ we set
\begin{equation}
\begin{array}{lcl}
\D(L_{\Omega}(t))&:=& \{u \in W^{2,p}(\R^d)\cap W^{1,p}_0(\Omega): \langle M(t)x, \mathrm D_xu(x)
\rangle \in L^p(\Omega) \},\\[0.2cm]
L_{\Omega}(t)u &:=& \OU_{\Omega}(t) u.
\end{array}
\end{equation}
Here the domain of $L_{\Omega}(t)$ depends on the time parameter $t$, however the subspace
\begin{equation*}
Y_{\Omega}:=\{u \in W^{2,p}(\Omega)\cap W^{1,p}_0(\Omega) : |x|\cdot \mathrm D_j u(x) \in L^p(\Omega) \; \mathrm{for} \; j=1,\ldots, d\}
\end{equation*}
is contained in $ \D(L_\Omega(t))$ for all $t\ge 0$ and is dense in $L^p(\Omega)$. It follows from \cite{Geissert/etal:2005} that in the autonomous case (i.e. for fixed $s\ge 0$) the operator $(L_{\Omega}(s),\D(L_{\Omega}(s))$ generates a strongly continuous semigroup on $L^p(\Omega)$. For more general second order elliptic operators with unbounded and time-independent coefficients in exterior domains we refer to \cite{Hieber/etal:2007}. Our main result is the existence of an evolution system in $L^p(\Omega)$, $1<p<\infty$, associated to the operators $L_\Omega(\cdot)$.
\begin{theorem}\label{prop:evolution_system_exterior}
Let $\Omega \subset \R^d$ be an exterior domain with $C^{1,1}$-boundary and $1<p<\infty$. Then there exists a unique evolution system $\{P_{\Omega}(t,s)\}_{(t,s)\in \Lambda}$ on $L^p(\Omega)$ with the following properties.
\begin{itemize}
\item[(a)] For $(t,s)\in \Lambda$, the operator $P_{\Omega}(t,s)$ maps $Y_{\Omega}$ into $Y_{\Omega}$.
\item[(b)] For every $f\in Y_{\Omega}$ and $s\ge 0$, the map $t\mapsto P_{\Omega}(t,s)f$ is differentiable in $(s,\infty)$ and
\begin{equation}\label{eq:derivative_evolution_exterior}
\frac{\partial}{\partial t} P_{\Omega}(t,s)f = L_{\Omega}(t)P_{\Omega}(t,s)f .
\end{equation}
\item[(c)] For every $f\in Y_{\Omega}$ and $t>0$, the map $s\mapsto P_{\Omega}(t,s)f$ is differentiable in $[0,t)$ and
\begin{equation}\label{eq:derivative_evolution_exterior_2}
\frac{\partial}{\partial s} P_{\Omega}(t,s)f = - P_{\Omega}(t,s)L_{\Omega}(s)f.
\end{equation}
\end{itemize}
\end{theorem}
As a direct consequence we obtain well-posedness of the abstract non-autonomous Cauchy problem (nACP) on the regularity space $Y_\Omega$.
\begin{corollary}
Let $\Omega$ be an exterior $C^{1,1}$-domain. Then the Cauchy problem (nACP) is well-posed on $Y_\Omega$.
\end{corollary}
In the following, we describe the construction of the evolution system $\{P_{\Omega}(t,s)\}_{(t,s)\in \Lambda}$ in detail. The general idea is to derive the result for exterior domains from the corresponding results in the case of $\R^d$ and bounded domains. For this purpose let $R>0$ be such that $K \subset B(R)$. We set $D:=\Omega \cap B(R+3)$. We denote by $\{P_{\R^d}(t,s)\}_{(t,s)\in \Lambda}$ the evolution system in $L^p(\R^d)$ and by $\{P_{D}(t,s)\}_{(t,s)\in \Lambda}$
the evolution system in $L^p(D)$ for the bounded domain $D$. Next we choose cut-off functions $\varphi,\eta \in C^{\infty}(\Omega)$ such that $0\leq \varphi, \eta \leq 1$ and
$$
\varphi(x) :=
\left\{
\begin{array}{cc}
  1,&|x| \geq R+2,  \\
  0, & |x| \leq R+1,
\end{array}
\right.
$$
and
$$
\eta(x) :=
\left\{
\begin{array}{cc}
  1,&|x| \leq R+2,  \\
  0, & |x| \geq R+\frac 5 2.
\end{array}
\right.
$$
For $f\in L^p(\Omega)$ we define $f_0 \in L^p(\R^d)$ and $f_D \in L^p(D)$, respectively, by
$$
f_0(x) :=
\left\{
\begin{array}{cc}
  f(x),&x\in\Omega,   \\
  0, & x\not\in \Omega,
\end{array}
\right. \qquad \mbox{and}\qquad f_D(x) = \eta(x) f(x).
$$
These definitions ensure that for every function $f\in \D(L_\Omega(t))$ we have $f_0\in \D(L_{\R^d}(t))$ and $f_D\in \D(L_D(t))$.
Now for $(t,s)\in \Lambda$ and $f\in L^p(\Omega)$ we set
\begin{equation}\label{eq:def_W(t,s)}
W(t,s)f = \varphi P_{\R^d}(t,s)f_0 + (1-\varphi) P_{D}(t,s)f_D.
\end{equation}
A short calculation yields
\begin{align*}
\mathrm D_x W(t,s)f & = \varphi \mathrm D_x P_{\R^d}(t,s)f_0 + (1-\varphi) \mathrm D_x P_{D}(t,s)f_D \\
&\qquad + \mathrm D_x \varphi \left(P_{\R^d}(t,s)f_0 - P_{D}(t,s)f_D\right),
\end{align*}
and
\begin{align*}
\mathrm D_x^2 W(t,s)f  & = \varphi \mathrm D_x^2 P_{\R^d}(t,s)f_0 + (1-\varphi) \mathrm D_x^2 P_{D}(t,s)f_D \\
& \qquad + 2 \,\left(\mathrm D_x \varphi\right)^* \cdot \left( \mathrm D_x P_{\R^d}(t,s)f_0 - \mathrm D_x P_{D}(t,s)f_D\right) \\
&\qquad + \mathrm D_x^2 \varphi \left(P_{\R^d}(t,s)f_0 - P_{D}(t,s)f_D\right).
\end{align*}
Thus, for $f\in Y_\Omega$, we obtain
\begin{equation}\label{eq:derivative_cutoff}
\left\{\begin{array}{rcll}
 \frac{\partial}{\partial t} W(t,s)f &=& L_\Omega(t)W(t,s)f - F(t,s)f, & (t,s)\in \Lambda ,\\[0.2cm]
 W(s,s)f&=&f, &
 \end{array}\right.
\end{equation}
with
\begin{align}\label{eq:error_terms}
F(t,s)f & = \mathrm{Tr}\left[Q(t)Q^*(t)\left(\mathrm D_x \varphi\right)^*\cdot \left( \mathrm D_x P_{\R^d}(t,s)f_0 - \mathrm D_x P_{D}(t,s)f_D\right)\right] \\
& \qquad + \OU_\Omega(t) \varphi\left(P_{\R^d}(t,s)f_0 - P_{D}(t,s)f_D\right).\notag
\end{align}
From the properties of the evolution systems $\{P_{\R^d}(t,s)\}_{(t,s)\in \Lambda}$ and $\{P_D(t,s)\}_{(t,s)\in \Lambda}$ it follows that the function $F(t,s)f$ in (\ref{eq:error_terms}) is well-defined for every $f\in L^p(\Omega)$ and $(t,s)\in \widetilde{\Lambda}$. Moreover, for every $f\in L^p(\Omega)$, $F(\cdot,\cdot)f$ is continuous in $\widetilde{\Lambda}$
with values in $L^p(\Omega)$. By using Proposition \ref{prop:Lp_Lq_estimates_whole} and Corollary \ref{prop:Lp_Lq_estimates_bounded} we obtain the estimate
\begin{equation}\label{eq:error-terms}
\|F(t,s)f\|_{p} \leq C \left(1+ (t-s)^{-\frac 1 2}\right) \|f\|_{p}, \qquad (t,s)\in \widetilde{\Lambda}_T,
\end{equation}
for any $T>0$ and a suitable constant $C:=C(T)>0$.

It is clear, that if an evolution system $\{P_{\Omega}(t,s)\}_{(t,s)\in \Lambda}$ exists on $L^p(\Omega)$ , then the solution $u(t)$ to the inhomogeneous equation (\ref{eq:derivative_cutoff}) is given by the variation of constant formula
$$
u(t) = P_{\Omega}(t,s)f - \int_s^t P_{\Omega}(t,r) F(r,s)f \d r.
$$
This consideration suggests to consider the integral equation
\begin{equation}\label{eq:integral}
P_\Omega(t,s)f=W(t,s)f+\int_s^tP_\Omega(t,r)F(r,s)f\,dr,\quad (t,s)\in \Lambda ,\,f\in L^p(\Omega).
\end{equation}
Let us state a lemma which will be very useful. Its proof is analogous to the proof in the case of one-parameter families (see \cite[Lemma 4.6]{Geissert/Heck/Hieber:2006}). But for the sake of completeness we give here the
details of the proof.
\begin{lemma}\label{lemma:iterated_convolution}
Let $X_1$ and $X_2$ be two Banach spaces, $T>0$ and let $R:\widetilde{\Lambda}_T \to \L(X_2,X_1)$ and $S:\widetilde{\Lambda}_T \to \L(X_2)$ be strongly continuous functions. Assume that
\begin{equation*}
\|R(t,s)\|_{\L(X_2,X_1)} \leq C_0(t-s)^{\alpha}, \quad \|S(t,s)\|_{\L(X_2)} \leq C_0(t-s)^{\beta}, \quad (t,s)\in \widetilde{\Lambda}_T,
\end{equation*}
holds for some $C_0:=C_0(T)>0$ and $\alpha, \beta >-1$. For $f\in X_2$ and $(t,s)\in \widetilde{\Lambda}_T$, set
$T_0(t,s)f:=R(t,s)f$ and
\begin{equation*}
T_n(t,s)f:= \int_s^t T_{n-1}(t,r)S(r,s)f\d s, \qquad n\in \N, \; (t,s)\in \widetilde{\Lambda}_T.
\end{equation*}
Then there exists a constant $C > 0$ such that
\begin{equation}\label{eq:convergence_series}
\sum_{n=0}^\infty \|T_n(t,s)f\|_{X_1} \leq C (t-s)^{\alpha}\|f\|_{X_2}, \qquad (t,s)\in \widetilde{\Lambda}_T.
\end{equation}
Moreover, if $\alpha \geq 0$, the convergence of the series in (\ref{eq:convergence_series}) is uniform on $\Lambda_T$.
\end{lemma}
\begin{proof}
For $f\in X_2$ and $(t,s)\in \widetilde{\Lambda}_T$ we have
$$\|T_1(t,s)f\|_{X_1}\le C_0^2\int_s^t(t-r)^\alpha (r-s)^\beta \,dr =C_0^2(t-s)^{\alpha +\beta +1}\B(\beta +1,\alpha +1)\|f\|_{X_2},$$
where $\B(\cdot ,\cdot)$ denotes the Beta function. So, by induction, we obtain
\begin{eqnarray*}
& & \|T_n(t,s)f\|_{X_1} \\
&\le & C_0^{n+1}(t-s)^{\alpha +n(\beta +1)}\B(\beta +1,\alpha +1)\cdots \B(\beta +1,\alpha +1+(n-1)(\beta +1))\|f\|_{X_2}\\
&=& C_0^{n+1}(t-s)^{\alpha +n(\beta +1)}\Gamma(\beta +1)^n\frac{\Gamma(\alpha +1)}{\Gamma(\alpha +1+n(\beta +1))}\|f\|_{X_2},\,\,n\in \N ,\,(t,s)\in \widetilde{\Lambda}_T,
\end{eqnarray*}
where $\Gamma$ denotes the Gamma function. Let us recall now the identity $\Gamma(x+1)=x\Gamma(x),\,x>-1$, and denotes by $[\cdot]$ the Gaussian brackets. Then, it follows that
$$\frac{\Gamma(\alpha +1)}{\Gamma(\alpha +1+n(\beta +1))}\le \frac{C_\alpha}{[n(\beta +1)]!},\quad n\in \N$$
for some $C_\alpha >0$.
Hence,
\begin{eqnarray*}
\|T_n(t,s)f\|_{X_1} &\le & C_\alpha C_0(t-s)^{\alpha }\Gamma(\beta +1)^nC_0^n\frac{(t-s)^{n(\beta +1)}}{[n(\beta +1)]!}\|f\|_{X_2}\\
&\le & C_\alpha C_0(t-s)^{\alpha }e^{t-s}\left(C_0\Gamma(\beta +1)\right)^n\frac{(t-s)^{[n(\beta +1)]}}{[n(\beta +1)]!}\|f\|_{X_2},\,
n\in \N ,\,(t,s)\in \widetilde{\Lambda}_T.
\end{eqnarray*}
Since
\begin{eqnarray*}
\sum_{n=0}^\infty \left(C_0\Gamma(\beta +1)\right)^n\frac{(t-s)^{[n(\beta +1)]}}{[n(\beta +1)]!} &\le & C_\beta e^{c_\beta(t-s)}\\
 &\le & C_\beta e^{c_\beta T}=:C_T,\quad (t,s)\in \Lambda_T
 \end{eqnarray*}
 for some constants $C_\beta ,\,c_\beta >0$, it follows that
 $$\sum_{n=0}^\infty \|T_n(t,s)f\|_{X_1}\le C_TC_0C_\alpha e^Tt-s)^{\alpha }\|f\|_{X_2},\quad (t,s)\in \widetilde{\Lambda}_T.$$
It is clear that if $\alpha \ge 0$ then the convergence of the above series is uniform on $\Lambda_T$.
\end{proof}

\begin{proof}[Proof of Theorem \ref{prop:evolution_system_exterior}.]
Let $T>0$. By using Proposition \ref{prop:Lp_Lq_estimates_whole} and Corollary \ref{prop:Lp_Lq_estimates_bounded} we have
\begin{equation*}
\|W(t,s)f\|_{p} \leq C\|f\|_p, \qquad \mbox{for}\; f \in L^p(\Omega), \; (t,s)\in \Lambda_T.
\end{equation*}
So, by (\ref{eq:error-terms}), we can apply Lemma \ref{lemma:iterated_convolution} with $R=W$, $S=F,\,\alpha =0,\,\beta =-\frac{1}{2}$ and $X_1=X_2=L^p(\Omega)$. Thus, for any $f\in L^p(\Omega)$,
the series $\sum_{k=0}^{\infty}P_k(t,s)f$ converges uniformly in $\Lambda_T$,
where $P_0(t,s)f=W(t,s)f$ and
\begin{equation}
P_{k+1}(t,s)f= \int_s^t P_k(t,r)F(r,s)f \d r,\quad (t,s)\in \Lambda_T,\,f\in L^p(\R^d).
\end{equation}
Since $T>0$ is arbitrary,
\begin{equation}\label{eq:representation_evolution_exterior}
P_\Omega (t,s):=\sum_{k=0}^{\infty}P_k(t,s),\quad(t,s)\in \Lambda
\end{equation}
is well-defined. It is easy to check that $P_\Omega(t,s)$ satisfies the integral equation \eqref{eq:integral}. Moreover, from the strong continuity of $W(\cdot ,\cdot)$ and (\ref{eq:error-terms}) we deduce inductively that $P_k(\cdot ,\cdot)$ is strongly continuous
and hence, by the uniform convergence of the series we get the strong continuity of $P_\Omega(\cdot ,\cdot)$.

In order to show that $\{P_{\Omega}(t,s)\}_{(t,s)\in \Lambda}$ leaves $Y_\Omega$ invariant, we consider the Banach space
$X_1:=\{f\in W_0^{1,p}(\Omega) : |x|\cdot \mathrm D_j f(x) \in L^p(\Omega)\; \mbox{for} \; j=1,\ldots,d\}$ endowed with the norm
$$\|f\|_{X_1}:=\|f\|_{1,p}+\||x|\cdot \mathrm D_xf\|_{p},\quad f\in X_1.$$
 Proposition \ref{prop:Lp_Lq_estimates_whole}, Corollary \ref{prop:Lp_Lq_estimates_bounded} and the last part of the proof of Proposition \ref{prop:evolution_whole} permit us to apply
 Lemma \ref{lemma:iterated_convolution} with $X_2=X_1,\,R=W,\,S=F,\alpha =0$ and $\beta =-\frac{1}{2}$. So,
 we obtain that $P_\Omega(t,s)f \in X_1$ for all $f\in X_1$ and $(t,s)\in \Lambda$. Moreover, by taking $X_1=W^{2,p}(\Omega),\,X_2=W^{1,p}(\Omega),\,R=W,\,S=F,\,\alpha=\beta=-\frac{1}{2}$ and applying Proposition \ref{prop:Lp_Lq_estimates_whole} and Corollary \ref{prop:Lp_Lq_estimates_bounded}, it follows, by Lemma \ref{lemma:iterated_convolution}, that
  $P_\Omega(t,s)f \in W^{2,p}(\Omega)$ for all $f\in W^{1,p}(\Omega)$ and $(t,s)\in \widetilde{\Lambda}$. This yields that $\{P_{\Omega}(t,s)\}_{(t,s)\in \Lambda}$ leaves $Y_\Omega$ invariant and
  \begin{align}\label{conv-in-Y_Omega}
  &\sum_{n=0}^\infty \left[\|P_k(t,s)f\|_{2,p}+\||x|\mathrm D_xP_k(t,s)f\|_p\right]\notag\\
  &\qquad\quad\quad <C_T(1+(t-s)^{-\frac{1}{2}})
   (\|f\|_{1,p}+\||x|\cdot \mathrm D_xf\|_{p}),\quad (t,s)\in \widetilde{\Lambda}_T ,\,f\in Y_\Omega.
  \end{align}

Let us now prove Equation (\ref{eq:derivative_evolution_exterior}). For $f\in Y_\Omega$ we compute
\begin{align*}
\frac{\partial}{\partial t} P_0(t,s)f & = L_\Omega(t) P_0(t,s) f - F(t,s)f\\
\frac{\partial}{\partial t} P_1(t,s)f & = L_\Omega(t) P_1(t,s) f + F(t,s)f - \int_s^t F(t,r)F(r,s)f \d r\\
\frac{\partial}{\partial t} P_2(t,s)f & = L_\Omega(t) P_2(t,s) f +  \int_s^t F(t,r)F(r,s)f \d r \\
& \qquad - \int_s^t \int_{r_1}^t F(t,r_2)F(r_2,r_1)F(r_1,s)f \d r_2 \d r_1.
\end{align*}
Inductively we see that
\begin{equation}\label{eq:derivative}
\frac{\partial}{\partial t} \sum_{k=0}^n P_k(t,s)f  = L_\Omega(t) \sum_{k=0}^n P_k(t,s)f - R_n(t,s)f
\end{equation}
holds for $n\in \N$, where
\begin{equation*}
R_n(t,s)f := \int_s^t \int_{r_1}^t \ldots \int_{r_{n-1}}^t F(t,r_n)F(r_n,r_{n-1})\ldots F(r_1,s)f\d r_n \ldots \d r_2 \d r_1.
\end{equation*}
Now, we estimate the norm of $R_n(t,s)f$. Estimate \eqref{eq:error_terms} yields
\begin{eqnarray*}
\|R_1(t,s)f\|_{p}&\leq & C^2 \int_s^t (t-r)^{-\frac 1 2}(r -s)^{-\frac 1 2} \d r \|f\|_p= C^2 \B(1/2,1/2)\|f\|_p,\\
\|R_2(t,s)f\|_{p} & \leq & C^3 \B(1/2,1/2)\int_s^t (r-s)^{-\frac{1}{2}}\d r \|f\|_p\\
&=& C^3 \B(1/2,1/2)\B(1/2,1)(t-s)^{\frac 1 2}\|f\|_p .
\end{eqnarray*}
Inductively, we see that
\begin{align}\label{eq:correction_terms}
\|R_n(t,s)\|_{p} &\leq C^{n+1}\B(1/2,1/2)\B(1/2,1)\ldots \B(1/2,n/2)(t-s)^{\frac{n-1}{2}}\|f\|_p\nonumber\\
& \leq \frac{C^{n+1}\Gamma(1/2)^n}{\left[\frac{n-1}{2}\right]!}(t-s)^{\frac{n-1}{2}}\|f\|_p
\end{align}
holds for $n\in \N$. Here the constant $C$ may change from line to line. From estimate (\ref{eq:correction_terms}) it follows that $\|R_n\|_{p} $ tends to zero as $n\rightarrow \infty$. So, by (\ref{conv-in-Y_Omega}) and
the closedness of $L_\Omega(t)$, we can conclude that
\begin{equation*}
\frac{\partial}{\partial t}P_\Omega(t,s)f  = L_\Omega(t) \sum_{k=0}^{\infty} P_k(t,s)f,\quad t>s,\,f\in Y_\Omega ,
\end{equation*}
 holds and this proves (\ref{eq:derivative_evolution_exterior}).

 Let us now show Equation (\ref{eq:derivative_evolution_exterior_2}). For $f\in Y_\Omega$ we have
 \begin{equation*}
L_D(s)(\eta f)= \eta L_\Omega(s) f + \mathrm{Tr}[Q(t)Q^*(t)(\mathrm D_x \eta)^* \cdot \mathrm D_x f] + (\OU_\Omega(s)\eta)f
\end{equation*}
holds. Thus,
\begin{align*}
W(t,s)L_\Omega (s)f&= \varphi P_{\R^d}(t,s)(L_\Omega(s)f)_0 + (1-\varphi) P_D(t,s)(L_\Omega (s)f)_D\\
& = \varphi P_{\R^d}(t,s)L_{\R^d}(s)f_0 + (1-\varphi) P_D(t,s)L_D (s)f_D - G(t,s)f,
\end{align*}
where
\begin{equation*}
G(t,s)f:= (1-\varphi)P_D(t,s)\left(\mathrm{Tr}[Q(t)Q^*(t)(\mathrm D_x \eta)^* \cdot \mathrm D_x f] + (\OU_\Omega(s)\eta)f\right)
\end{equation*}
and $f\in Y_\Omega$.
This yields
\begin{equation*}
\frac {\partial}{\partial s} W(t,s)f = - W(t,s)L_\Omega (s)f - G(t,s)f
\end{equation*}
for $(t,s)\in \Lambda$ and $f\in Y_\Omega$.

Now, let $T>0$ be arbitrary but fixed. Then, from the definition of $G$ and Corollary \ref{prop:Lp_Lq_estimates_bounded}, it follows that we can apply Lemma \ref{lemma:iterated_convolution} with $X_1=X_2=W^{1,p}(\Omega),\,R=S=G$ and $\alpha =\beta =-\frac{1}{2}$. So, the series
 \begin{equation*}
T(t,s)f := \sum_{k=0}^{\infty} T_k (t,s)f, \qquad (t,s)\in \widetilde{\Lambda}_T,
\end{equation*}
is well-defined and
\begin{equation}\label{eq-22}
\|T(t,s)f\|_{1,p}\le C(t-s)^{-\frac{1}{2}}\|f\|_{1,p},\quad (t,s)\in \widetilde{\Lambda}_T,
\end{equation}
 for $f\in W^{1,p}(\Omega)$. On the other hand, $T(\cdot,\cdot)$ satisfies the integral equation
\begin{equation}\label{equ:integral}
T(t,s)f = G(t,s)f + \int_s^t T(t,r)G(r,s)f \d r, \qquad (t,s)\in \Lambda_T,\,f\in W^{1,p}(\Omega).
\end{equation}
In particular $T(t,\cdot)f$ is continuous on $[0,t]$ with respect to the $L^p$-norm for any $f\in W^{1,p}(\Omega)$ and
$t\ge 0$. Now, for $f\in L^p(\Omega)$ and $(t,s)\in \Lambda_T$ we set
\begin{equation*}
S(t,s)f := W(t,s)f+\int_s^t T(t,r)W(r,s)f \d r.
\end{equation*}
It follows from the continuity of $T(t,\cdot)W(\cdot ,s)f$ on $[s,t]$, Proposition \ref{prop:Lp_Lq_estimates_whole} and Corollary \ref{prop:Lp_Lq_estimates_bounded} that the above integral is well-defined for any $f\in L^p(\Omega)$.
Computing the derivative with respect to $s$ yields
\begin{align*}
\frac{\partial}{\partial s} S(t,s)f &= - W(t,s)L_\Omega(s)f - G(t,s)f+ T(t,s)f - \int_s^t T(t,r) W(r,s)L_\Omega(s) f \d r \\
& \qquad \quad - \int_s^t T(t,r)G(r,s)f \d r  \\[0.2cm]
& = - S(t,s)L_\Omega(s)f ,
\end{align*}
for $f\in Y_\Omega$, due to (\ref{equ:integral}). From this equality together with (\ref{eq:derivative_evolution_exterior}) and since
$P_\Omega(t,s)Y_\Omega \subset Y_\Omega,\,(t,s)\in \Lambda$, we can conclude that
\begin{align*}
\frac{\partial}{\partial r} (S(t,r)P_\Omega(r,s)f)= 0
\end{align*}
holds for all $f\in Y_\Omega$. This yields that for $f\in Y_\Omega$, the function $S(t,r)P_\Omega(r,s)f$ is constant on $\Lambda_T $ and thus, by the density of $Y_\Omega$ in $L^p(\Omega)$ and by the fact that $T>0$ was arbitrary, it follows that $S(t,s)f=P_\Omega(t,s)f$ holds for all $f\in L^p(\Omega)$ and all $(t,s)\in \Lambda$. This proves  (\ref{eq:derivative_evolution_exterior_2}).

Let us now show the uniqueness of the solution $P_\Omega(t,s)f$ of (nACP) for initial value $f \in Y_\Omega$. For this purpose we assume that there exists another solution $t
\mapsto u(t) \in Y_\Omega$. Since $u(r) \in Y_\Omega$ for all $r\in[s,\infty)$ it follows from equality (\ref{eq:derivative_evolution_exterior_2}) that the map $r \mapsto P_{\Omega}(t,r)u(r)$ is differentiable  for $0\leq s< r < t $ and
\begin{equation*}
\frac{\partial}{\partial r}\left( P_{\Omega}(t,r)u(r)\right) = -P_{\Omega}(t,r)L_\Omega(r)u(r) + P_{\Omega}(t,r)L_\Omega(r)u(r) = 0.
\end{equation*}
Therefore $P_{\Omega}(t,r)u(r)$ is constant on $0\leq s< r < t $. Thus, by letting $r\to s$ and $r\to t$ we obtain $P_\Omega(t,s)f = u(t)$. The uniqueness now directly implies that the law of evolution (Property (i) of Definition \ref{def:evolution_system}) holds.
\end{proof}

To conclude this section we prove $L^p$-$L^q$ smoothing properties of the evolution system $\{P_{\Omega}(t,s)\}_{(t,s)\in\Lambda}$ and $L^p$-estimates for its spatial derivatives. The following estimates follow basically directly via the representation \eqref{eq:representation_evolution_exterior} from Lemma \ref{lemma:iterated_convolution}, Proposition \ref{prop:Lp_Lq_estimates_whole}  and Corollary  \ref{prop:Lp_Lq_estimates_bounded}.
\begin{prop}
Let $T>0$, $1<p<\infty$ and $p\leq q < \infty$. Then there exists a constant $C:=C(T)>0$ such that
\begin{enumerate}
\item[(i)] $\|P_{\Omega}(t,s)f\|_{q} \leq C(t-s)^{-\frac{d}{2}\left(\frac 1 p - \frac 1 q\right)}\|f\|_{p}$,
\item[(ii)] $\|\mathrm D_x P_\Omega(t,s)f\|_p\le C(t-s)^{-\frac{1}{2}}\|f\|_p$
\end{enumerate}
for $(t,s)\in \widetilde{\Lambda}_T$ and $f\in L^p(\Omega)$. Moreover, for $1<p<q<\infty$ and $f\in L^p(\Omega)$
$$\lim_{t\to s}\left[\|(t-s)^{\frac{d}{2}\left(\frac 1 p - \frac 1 q\right)}P_\Omega(t,s)f\|_q+\|(t-s)^{\frac{1}{2}}\mathrm D_x P_\Omega(t,s)f\|_p\right]=0.$$
\end{prop}
\begin{proof}
To obtain (i) we apply Lemma \ref{lemma:iterated_convolution} with $X_1=L^q(\Omega),\,X_2=L^p(\Omega),\,R=W,\,S=F,\,\alpha=
-\frac{d}{2}\left(\frac 1 p - \frac 1 q\right),\,\beta =-\frac{1}{2}$,
Proposition \ref{prop:Lp_Lq_estimates_whole}  and Corollary  \ref{prop:Lp_Lq_estimates_bounded} in the case where $q\ge p\ge \frac{d}{2}$. By iteration (i) holds also for $1<p<\frac{d}{2}$.

The second assertion follows by applying Lemma \ref{lemma:iterated_convolution} with $X_1=W^{1,p}(\Omega),\,X_2=L^p(\Omega),\,R=W,\,S=F,\,\alpha =\beta =-\frac{1}{2}$, Proposition \ref{prop:Lp_Lq_estimates_whole}  and Corollary  \ref{prop:Lp_Lq_estimates_bounded}. Finally, the last assertion can be obtained as in \cite[Proposition 3.4]{Hieber-Sawada}.
\end{proof}
%
%
%\begin{thebibliography}{99}
\bibliographystyle{alpha}
%\bibliography{references}

%
%
%
%\end{thebibliography}
%
\end{document}